\documentclass[10pt]{article}

\usepackage{a4wide}
\usepackage{amssymb}
\usepackage{amsfonts}
\usepackage{amsmath}
\input xy
\xyoption{arrow} \xyoption{matrix}

\date{}

\newtheorem{proposition}{Proposition}[section]
\newtheorem{theorem}[proposition]{Theorem}
\newtheorem{lemma}[proposition]{Lemma}

\newtheorem{corollary}[proposition]{Corollary}

\def\der{\partial }

\def\nFM0{{\nu }_{F,M_0}}
\def\nFN0{{\nu }_{F,N_0}}
\def\nGN0{{\nu }_{G,N_0}}

\def\N0{ {\bf N}_0 }

\def\t{\otimes}
\def\g{\gamma}

\def\ra{\rightarrow}

\def\Xpm{X^{\pm }}

\def\s{\sigma}
\def\Z{\mathbb{Z}}

\def\l1{{\lambda}_1}

\def\a{\alpha}
\def\a0{ {\alpha }_0}
\def\a1{ {\alpha }_1}

\def\l{\lambda}

\def\nFGM0{{\nu }_{F,G,M_0}}

\def\lc{{\rm lc}}
\def\nFN0{{\nu}_{F,N_0}}

\def\sm{{\sigma}^m}

\def\sm1{{\sigma}^{-1}}

\def\smtp1{{\sigma}^{-t+1}}

\def\S1{S^{-1}}

\def\Xpm1{X^{\pm 1}_1}

\def\sPM1{{\sigma }^{\pm 1}}
\def\sMP1{{\sigma }^{\mp 1 }}

\def\b{\beta}
\def\d{\delta}

\def\di{{\rm d.ind}}

\def\L{\Lambda}

\def\CA{{\cal A}}

\def\CD{{\cal D}}

\def\Ytm1{Y^{t-1}}
\def\Yim1{Y^{i-1}}

\def\CF{{\cal F}}

\def\CH{{\cal H}}

\def\CZ{{\cal Z}}
\def\supp{{\rm supp }}

\def\Der{{\rm Der }}
\def\ad{{\rm ad }}

\def\pad{{\rm pad }}

\def\SL2Z{ {\rm SL}_2({\bf Z}) }

\def\CZ{ {\cal Z}}

\def\Gp1{ G^{1 , 1 } }
\def\P11{ P^{-1 , 1 } }
\def\Pp1{ P^{1 , 1 } }

\def\nCLsr{{}^\nu\kern-2pt {\cal L}^{\sigma , \rho  }}
\def\nP{{}^\nu \kern-2pt P}
\def\nL{{}^\nu\kern-2pt L}
\def\nLL{{}^\nu\kern-2pt \Lambda}
\def\nPsr{{}^\nu\kern-2pt P^{\sigma , \rho  }}
\def\nLsr{{}^\nu\kern-2pt L^{\sigma , \rho  }}
\def\nuCL{{}^\nu\kern-2pt  {\cal L}}
\def\nCLsr{{}^\nu\kern-2pt {\cal L}^{\sigma , \rho  }}
\def\nCL1m{{}^\nu\kern-2pt {\cal L}^{-1 , 1  }}
\def\x1nu{x^\frac{1}{\nu}}
\def\xm1nu{x^{-\frac{1}{\nu}}}

\def\ob{\overline{b}}

\def\ra{\rightarrow }

\def\CB{{\cal B}}

\def\CH{ {\cal H}}

\def\nAM0{{\nu }_{{\cal A},M_0}}
\def\nAN0{{\nu }_{{\cal A},N_0}}

\def\Der{ {\rm Der }}

\def\ad{ {\rm ad }}

\def\ga{\mathfrak{a}}
\def\gb{\mathfrak{b}}

\def\SL{{\rm SL}}

\def\di!{\frac{\der^i}{i!}}
\def\dik!{\frac{\der^k_i}{k!}}

\def\gl{\mathfrak{l}}

\def\id{{\rm id}}

\def\N{\mathbb{N}}

\def\0{\overline{0}}
\def\1{\overline{1}}

\def\Ln1{\L_{n,\overline{1}}}

\def\oa{\overline{a}}

\def\IDer{{\rm IDer}}

\def\a1{a_{\overline{1}}}

\def\S{\Sigma}

\def\sign{{\rm sign}}
\def\vn1{\overrightarrow{n-1}}

\def\gl{{\rm gl}}
\def\sl{{\rm sl}}

\def\bu{\overline{u}}
\def\PZ{{\rm PZ}}

\def\mJ{\mathbb{J}}
\def\mI{\mathbb{I}}

\def\ann{{\rm ann}}

\def\mT{\mathbb{T}}

\def\K1{{\rm K}_1}

\def\hmI1{\widehat{\mI_1}}
\def\tmI1{\widetilde{\mI_1}}
\def\tmJ1{\widetilde{\mJ_1}}
\def\hB1{\widehat{B_1}}
\def\hCB1{\widehat{\CB_1}}

\def\ga{\mathfrak{a}}

\def\sl2{\mathfrak{sl}_2}
\def\gl2{\mathfrak{gl}_2}

\def\b1{\overline{1}}

\def\PDer{{\rm PDer}}
\def\IDer{{\rm IDer}}
\def\PIDer{{\rm PIDer}}
\def\AutPois{{\rm Aut}_{\rm Pois}}
\def\lt{{\rm lt}}

\begin{document}

\author{V. V. \  Bavula 
}

\title{The generalized Weyl  Poisson algebras and their Poisson simplicity criterion}

\maketitle

\begin{abstract}

A new class of Poisson algebras, the class of   {\em  generalized Weyl Poisson algebras}, is introduced. It can be seen as Poisson algebra analogue of generalized Weyl algebras or as  giving a Poisson structure to (certain) generalized Weyl algebras. A Poisson simplicity criterion is given for  generalized Weyl Poisson algebras and explicit descriptions of  the Poisson centre and the absolute Poisson centre are obtained. Many examples are considered.

$\noindent $



$${\bf Contents}$$
\begin{enumerate}
\item Introduction.
\item The generalized Weyl  Poisson algebras.
\item Poisson simplicity criterion for  generalized Weyl  Poisson algebras
\end{enumerate}
\end{abstract}


\section{Introduction}\label{INTR}

$111$

In this paper, 
 $K$ is a field, algebra means a $K$-algebra (if it is not stated otherwise) and $K^*=K\backslash \{ 0\}$.



{\bf Generalized Weyl algebras, \cite{Bav-GWA-FA-91, Bav-SimGWA-1992, Bav-GWArep}.} Let $D$ be a ring, $\sigma=(\sigma_1,...,\sigma_n)$ be  an $n$-tuple of
commuting automorphisms of $D$,  $a=(a_1,...,a_n)$ be  an $n$-tuple  of elements of  the centre
$Z(D)$  of $D$ such that $\sigma_i(a_j)=a_j$ for all $i\neq j$. The {\em generalized Weyl algebra} $A=D[X, Y; \sigma,a]$ (briefly GWA) of rank  $n$  is  a  ring  generated
by $D$  and    $2n$ indeterminates $X_1,...,X_n, Y_1,...,Y_n$
subject to the defining relations:
$$Y_iX_i=a_i,\;\; X_iY_i=\sigma_i(a_i),\;\; X_id=\sigma_i(d)X_i,\;\; Y_id=\sigma_i^{-1}(d)Y_i\;\; (d \in D),$$
$$[X_i,X_j]=[X_i,Y_j]=[Y_i,Y_j]=0, \;\; {\rm for \; all}\;\; i\neq j,$$
where $[x, y]=xy-yx$. We say that  $a$  and $\sigma $ are the  sets  of
{\it defining } elements and automorphisms of the GWA $A$, respectively.

 The $n$'th {\em Weyl algebra},
 $A_n=A_n(K)$ over a field (a ring) $K$ is  an associative
 $K$-algebra generated by  $2n$ elements
 $X_1, ..., X_n,Y_1,..., Y_n$, subject to the relations:
$$[Y_i,X_i]=\d_{ij}\;\; {\rm and}\;\;  [X_i,X_j]=[Y_i,Y_j]=0\;\; {\rm for\,all}\; i,j, $$
where $\d_{ij}$ is the Kronecker delta function.
The Weyl algebra $A_n$ is a generalized Weyl algebra
 $A=D[X, Y; \s ;a]$ of rank $n$ where
$D=K[H_1,...,H_n]$ is a polynomial ring   in $n$ variables with
 coefficients in $K$, $\s = (\s_1, \ldots , \s_n)$ where $\s_i(H_j)=H_j-\delta_{ij}$ and
 $a=(H_1, \ldots , H_n)$.  The map
$$A_n\ra A,\;\; X_i\mapsto  X_i,\;\; Y_i \mapsto  Y_i,\;\;  i=1,\ldots ,n,$$
is an algebra  isomorphism (notice that $Y_iX_i\mapsto H_i$).

It is an experimental fact that many quantum algebras of small Gelfand-Kirillov dimension are GWAs (eg, $U({\rm sl}_2)$, $U_q({\rm sl}_2)$, the quantum Weyl algebra, the quantum plane, the Heisenberg algebra and its quantum analogues, the quantum sphere,  and many others).

The GWA-construction turns out to be a useful one. Using it for large classes of algebras (including the mention ones above) all the simple modules were classified, explicit formulae were found for the  global
and Krull dimensions, their elements were classified in the sense of Dixmier, \cite{Dix}, etc.

{\bf The  generalized Weyl Poisson algebra $D[X,Y; a, \der \}$.} Our aim is to introduce a Poisson algebra analogue of generalized Weyl algebras. Let $A$ be a Poisson algebra with   Poisson bracket $\{\cdot, \cdot\}$ , $\CZ (A):=\{ a\in A\, | \, ax=xa, \; \{a, x\}=0$ for all $x\in A\}$ be its {\em absolute centre} and $\PDer_K(A)$ be the set of derivations of the Poisson algebra $A$ (see Section \ref{PoisGWAs} for details).

 {\it Definition.} Let $D$ be a Poisson algebra (not necessarily commutative as an associative algebra), $\der = (\der_1, \ldots , \der_n)\in \PDer_K(D)^n$ be an $n$-tuple of commuting derivations of the Poisson algebra $D$, $a=(a_1, \ldots , a_n)\in \CZ (D)^n$ 
  be such that $\der_i (a_j)=0$ for all $i\neq j$. The generalized Weyl algebra $$A= D[ X, Y; (\id_D, \ldots , \id_D), a]=
 D[X_1, \ldots , X_n , Y_1, \ldots , Y_n]/(X_1Y_1-a_1, \ldots , X_nY_n-a_n)$$ admits a Poisson structure which is an extension of the Poisson structure on $D$ and is given by the rule: For all $i,j=1, \ldots , n$ and $d\in D$,

\begin{equation}\label{PGWAR1}
\{ Y_i, d\}=\der_i(d)Y_i, \;\; \{ X_i, d\}=-\der_i(d)X_i \;\; {\rm and}\;\; \{ Y_i, X_i\} = \der_i (a_i),
\end{equation}
\begin{equation}\label{PGWAR2}
\{ X_i, X_j\}=\{ X_i, Y_j\}=\{ Y_i, Y_j\} =0\;\; {\rm for\; all}\;\; i\neq j.
\end{equation}
The Poisson algebra is denoted by $A =D[ X, Y; a, \der \}$ and is called the {\bf generalized Weyl Poisson algebra} of rank $n$ (or GWPA, for short)  where $X=(X_1, \ldots , X_n)$ and $Y= (Y_1, \ldots , Y_n)$.

Existence of  generalized Weyl Poisson  algebras is proven in  Section \ref{PoisGWAs} (Lemma \ref{a10Apr16}). The key idea of the proof is to introduce another class of Poisson algebras, elements of which is denoted by  $D[ X,Y; \der , \alpha ]$ (see  Section \ref{PoisGWAs}),  for which existence problem has easy solution and then to show that each GWPA is a factor algebra of some $D[ X,Y; \der , \alpha ]$. The Poisson algebras $D[ X,Y; \der , \alpha ]$ turn out to be also GWPAs (Proposition \ref{b10Apr16}).

{\bf Poisson simplicity criterion for  generalized Weyl Poisson   algebras.} A Poisson algebra is a {\em simple} Poisson algebra if the ideals $0$ and $A$ of the associative algebra $A$ are the only ideals $I$ such that $\{ A, I\} \subseteq I$. The ideal $I$ is called a {\em Poisson ideal} of the Poisson algebra $A$. An ideal $I$ of the ring $D$ is called $\der$-{\em invariant}, where $\der = (\der_1, \ldots , \der_n)\in \PDer_K(D)^n$, if $\der_i(I)\subseteq I$ for all $i=1, \ldots , n$. The set $D^\der := \{ d\in D\, | \, \der_1(d)=0, \ldots , \der_n(d)=0\}$ is called the {\em ring of $\der$-constants} of $D$.

  In Section \ref{SIMCRPGWA}, a proof is given  of the following Poisson simplicity criterion for    generalized Weyl Poisson  algebras, see Proposition \ref{a13Apr16} for the notation.

\begin{theorem}\label{10Apr16}
Let $A=D[X,Y;a, \der\}$ be a GWPA of rank $n$. 
Then the  Poisson algebra $A$ is a simple Poisson algebra iff
\begin{enumerate}
\item the Poisson algebra  $D$ has no proper $\der$-invariant Poisson ideals,
\item for all $i=1, \ldots , n$, $Da_i+D\der_i(a_i) = D$, and
\item the algebra $\CZ (A)$ is a field, i.e. char$(K)=0$, $\CZ (D)^\der$ is a field and $D_{[\alpha ]}=0$ for all $\alpha \in \Z^n\backslash \{ 0\} $ (see the proposition below).
\end{enumerate}
\end{theorem}

As a first step in the proof of Theorem \ref{10Apr16}, the following field criterion for the absolute centre $\CZ (A)$ of a GWPA $A=D[X,Y;a, \der\}$  of rank $n$ is proven (in Section \ref{SIMCRPGWA}).

\begin{proposition}\label{bb13Apr16}
Let $A= D[X,Y; a, \der\}$ be a GWPA of rank $n$. Then $\CZ (A)$ is a field iff char$(K)=0$, $\CZ (D)^\der$ is a field and  $D_{[\alpha ]}=0$ for all $\alpha \in \Z^n\backslash \{ 0\} $.
\end{proposition}

An explicit descriptions of the Poisson centre and absolute centre are obtained (Proposition \ref{a13Apr16}).
Many examples are considered. We show that many classical Poisson algebras are GWPAs.

At the end of Section \ref{PoisGWAs}, we show  that GWPAs appear as associated graded Poisson algebras of certain GWAs (Proposition \ref{c10Apr16}). This is a sort of quantization procedure.

At the end of Section \ref{SIMCRPGWA}, examples of simple GWPAs (as Poisson algebras) are considered (Corollary \ref{a19Apr16}).  This family of simple Poisson algebras includes, as a particular case,
the {\em classical Poisson polynomial algebras} $P_{2n} = K[X_1, \ldots , X_n , Y_1, \ldots , Y_n]$ ($\{ Y_i, X_j\} =\d_{ij}$ and $\{X_i, X_j\} = \{ X_i, Y_j\} = \{ Y_i, Y_j\} =0$ for all $i\neq j$).


\section{The generalized Weyl  Poisson algebras}\label{PoisGWAs}

In this section, two new classes of Poisson algebras  are introduced and prove their existence. One of them is the class of   generalized Weyl Poisson  algebras (GWPAs). Examples are considered. At the end of the section, it is shown that some GWPAs are obtained from GWAs by a sort of quantization procedure (Proposition \ref{c10Apr16}).

{\bf Poisson algebras.} An associative (not necessarily commutative) algebra $D$ is called a {\em Poisson algebra} if it is a Lie algebra $(D, \{ \cdot, \cdot \})$ such that $\{ a, xy\}= \{ a, x\}y+x\{ a, y\}$ for all elements $a,x,y\in D$.

 For a $K$-algebra  $D$,  let $\Der_K(D)$ be the set of its $K$-derivations. If, in addition, $(D, \{ \cdot , \cdot \} )$ is a Poisson algebra then 
 $$\PDer_K(D):=\{ \d \in \Der_K(D)\, | \, \d (\{ a,b\})= \{ \d (a),b\}+\{ a,\d (b)\} \;\; {\rm  for\; all}\;\; a,b\in D\}$$
 is the set of derivations of the Poisson algebra $D$. The vector space $\Der_K(D)$ is a Lie algebra, where $[\d , \der ]:= \d \der -\der \d$, and $\PDer_K(D)$ is a Lie subalgebra of $\Der_k(D)$. The set of {\em inner derivations} 
 $$\IDer_K(D):= \{ \ad_a\, | \, a\in D\}\;\; {\rm  (where} \;\;  \ad_a (b) := [a,b]:=ab-ba)$$ 
 is an ideal of the Lie algebra $\Der_K(D)$ (since $[\d ,\ad_a]=\ad_{\d (a)}$ for all $\d\in \Der_K(D)$ and $a\in D$). Similarly, the set of {\em inner derivations of the Poisson algebra} $D$,  
 $$\PIDer_K(D):= \{ \pad_a\, | \, a\in D\}\;\; {\rm  (where}\;\;  \pad_a (b) := \{ a,b\} )$$ 
 is an ideal of the Lie algebra $\PDer_K(D)$ (since $[\d ,\pad_a]=\pad_{\d (a)}$ for all $\d\in \PDer_K(D)$ and $a\in D$). By the very definition, the Poisson algebra $D$ is a Lie algebra with respect to the bracket $\{ \cdot , \cdot \}$. The map $ D\ra \IDer_K(D)$, $a\mapsto \pad_a$, is an epimorphism of Lie algebras with kernel 
 $$\PZ (D):= \{ a\in D\, | \, \{ a, D\} =0\}$$ which is called the {\em centre} of the Poison algebra (or the {\em Poisson centre} of $D$). So, the Poisson  structure of the algebra $D$ induces the `multiplicative structure' on the Lie algebra $\PIDer_K(D)$, i.e. $\pad_{ab}(\cdot)= \pad_{a}(\cdot)\, b+a\, \pad_{b}(\cdot)$.

 Notice that the {\em centre} $Z(D) :=\{ z\in D\, | \, zd=dz$ for all $d\in D\}$ of any associative algebra $D$ is invariant under the action of $\Der_K(D)$: Let $z\in \PZ(D)$, $d\in D$ and $\der \in \Der_K(D)$; then applying the derivation $\der$  to the equality $zd=dz$ we obtain the equality $ \der (z) d= d\der (z)$, i.e. $\der (z) \in Z(D)$. Similarly,  the Poisson  centre $\PZ(D)$  is invariant under the action of $\PDer_K(D)$: Let $z\in \PZ(D)$, $d\in D$ and $\der \in \PDer_K(D)$; then applying the derivation  $\der$  to the equality $\{ z,d\} =0$ we obtain the equality  $\{ \der (z) ,d\}=0 $, i.e. $\der (z) \in \PZ(D)$. For a Poisson algebra $D$, the intersection 
 $$\CZ (D) := Z(D)\cap \PZ (D)$$ is called the {\em absolute centre} of the Poisson algebra $D$. {\em The absolute centre is invariant under the action of} $\PDer_K(D)$.

Let $D$ be a Poisson algebra where the associative algebra $D$ is not necessarily commutative. Let $\der = (\der_1, \ldots , \der_n)\in \PDer_K(D)^n$ be an $n$-tuple of commuting derivations of the Poisson algebra $D$ and $X= (X_1, \ldots , X_n)$ be an $n$-tuple of commuting variables. The polynomial algebra $D[X]=D[X_1, \ldots , X_n]$ with coefficients from $D$ admits a Poisson structure which is an extension of the Poisson structure on $D$ given by the rule
\begin{equation}\label{DXPA}
\{ X_i, X_j\}=0\;\; {\rm and}\;\; \{ X_i, d\} = \der_i(d) X_i\;\; {\rm for\;} 1\leq i,j\leq n\;\; {\rm and}\;\; d\in D.
\end{equation}
The Poisson algebra $D[X]$ is denoted by $D[ X; \der]$ and is called the {\em Poisson Ore extension} of $D$ of rank $n$. When $n=1$,  a more general construction appeared in *** Cho and Oh  ****

 Let $G$ be a monoid. Suppose that the associative algebra $D= \oplus_{g\in G} D_g$ is a $G$-{\em graded algebra}  ($D_gD_h\subseteq D_{gh}$ for all $g,h\in G$). If, in addition, $D$ is a Poisson algebra and  $\{ D_g, D_h\}\subseteq D_{gh}$ for all $g,h\in G$ then we say that the Poisson algebra $D$ is a $G$-{\em graded Poisson algebra}.

{\bf The Poisson algebra $D[X,Y; \der, \alpha ]$.} Now, we introduce a class of Poisson algebras which is used in the proof of existence of GWPAs (Lemma \ref{a10Apr16}).

 {\it Definition.} Let $D$ be a Poisson algebra, $\der = (\der_1, \ldots , \der_n)\in \PDer_K(D)^n$ be an $n$-tuple of commuting derivations of the Poisson algebra $D$ and $\alpha = (\alpha_1, \ldots , \alpha_n)\in \PZ (D)^n$. Then the polynomial algebra $D[X,Y]=D[X_1, \ldots , X_n, Y_1, \ldots , Y_n]$   with coefficients in $D$ admits a Poisson structure which is an extension of the Poisson structure on $D$ given by the rule: For all $i,j=1, \ldots , n$ and $d\in D$,

\begin{equation}\label{DXYPA1}
\{ Y_i, d\}=\der_i(d)Y_i, \;\; \{ X_i, d\}=-\der_i(d)X_i \;\; {\rm and}\;\; \{ Y_i, X_i\} = \alpha_i,
\end{equation}
\begin{equation}\label{DXYPA2}
\{ X_i, X_j\}=\{ X_i, Y_j\}=\{ Y_i, Y_j\} =0\;\; {\rm for\; all}\;\; i\neq j.
\end{equation}
The Poisson algebra $D[X,Y]$ is denoted by $\CA =D[ X, Y; \der , \alpha]$ where $X=(X_1, \ldots , X_n)$ and $Y= (Y_1, \ldots , Y_n)$.

Let us show that the Poisson structure on the polynomial algebra  $D[X,Y]$ is well-defined. Let $n=1$. The Poisson algebra $D[X_1, Y_1; \der_1, \alpha_1]$ is an extension of the Poisson  Ore extension $D[ X_1; -\der_1]$ by adding a commuting variable $Y_1$ where the Poisson structure on  the algebra $D[X_1][Y_1]$ is given by the rule $$\{ Y_1, d\} = \der_1(d) Y_1\;\; {\rm and}\;\;   \{ Y_1, X_1\} = \alpha_1.$$ The Poisson structure on  the algebra $D[X_1][Y_1]$ is well-defined as $\{ Y_1, \cdot\}$ respects the relation $\{ X_1, d\} = -\der_1(d) X_1$ for all $d\in D$:
$$ \{ Y_1, \{ X_1, d\}\}= \{ \alpha_1, d\} +\{ X_1, \der_1 (d) Y_1\}= 0 -\der_1^2(d)X_1Y_1-\der_1(d) \alpha_1=-\{ Y_1, \der_1(d) X_1\}. $$

For $n\geq 1$, the Poisson algebra
\begin{equation}\label{DXYPA3}
D[ X, Y; \der , \alpha]=D[ X_1, Y_1; \der_1 , \alpha_1]\cdots [ X_n, Y_n; \der_n , \alpha_n].
\end{equation}
is an iteration of this construction $n$ times.

Existence of the construction of   generalized Weyl Poisson  algebra follows from the next lemma.

\begin{lemma}\label{a10Apr16}
We keep the assumptions of the Definition of GWPA $A=D[X,Y; a, \der\}$. Let $\CA = D[X,Y; \der , \der (a)]$ where $\der (a) = (\der_1(a_1), \ldots , \der_n(a_n))$. Then $X_1Y_1-a_1, \ldots , X_nY_n-a_n\in \CZ (\CA )$ and the  generalized Weyl Poisson algebra
$A =D[ X, Y; a, \der \}$ is a factor algebra of the Poisson algebra $\CA$, $$A\simeq \CA /(X_1Y_1-a_1, \ldots , X_nY_n-a_n).$$
\end{lemma}

{\it Proof}. By the very definition,  the element $Z_i=X_iY_i-a_i\in Z(A)$: For all $i,j$ such that $i\neq j$, $\{ X_j, Z_i\} = \der_j(a_i) X_j=0$ and $\{ Y_j, Z_i\} = -\der_j(a_i) Y_j=0$ (since $\der_j(a_i)=0$ for all $i\neq j$). For all $d\in D$,
\begin{eqnarray*}
\{ Z_i, d\}  &=& \{ X_i, d\} Y_i+X_i\{ d, Y_i\} = -\der_i(d) X_iY_i+X_i\der (d) Y_i=0,  \\
\{ X_i, Z_i\}  &=& X_i(- \der_i(a_i))+\der_i(a_i) X_i=0, \\
\{ Y_i, Z_i\}  &=& \der_i(a_i)Y_i-\der_i(a_i) Y_i=0.
\end{eqnarray*}
Therefore,  $Z_i\in \PZ (\CA )$. Now, the lemma is obvious.  $\Box $

$\noindent $

The GWPA of rank $n$,
\begin{equation}\label{AZngr}
A:= D[X,Y; a, \der \} =\bigoplus_{\alpha\in \Z^n}A_\alpha ,
\end{equation}
is a $\Z^n$-graded Poisson algebra where $A_\alpha = Dv_\alpha$,  $v_\alpha =\prod_{i=1}^n v_{\alpha_i}(i)$ and $v_j(i)=\begin{cases}
X_i^j& \text{if }j>0,\\
1& \text{if }j=1,\\
Y_i^{|j|}& \text{if }j<0.\\
\end{cases}$

So, $A_\alpha A_\beta \subseteq A_{\alpha +\beta}$ and  $\{ A_\alpha , A_\beta \}\subseteq A_{\alpha +\beta}$ for all elements $\alpha , \beta \in \Z^n$.

{\bf The isomorphisms $s_I$ where $I\subseteq \{ 1, \ldots , n\}$ of  GWPAs of rank $n$.} Let $A=D[X_1,Y_1; a_1, \der_1\}$ be a GWPA of rank 1. Clearly, $A\simeq D[Y_1,X_1; a_1, -\der_1\}$, i.e.  the $D$-homomorphism of Poisson algebras
\begin{equation}\label{s1iso}
s_1: A=D[X_1,Y_1; a_1, \der_1 \}\ra D[Y_1,X_1; a_1, -\der_1\}, \;\; X_1\mapsto Y_1, \;\; Y_1\mapsto X_1, \;\; d\mapsto d\;\; (d\in D),
\end{equation}
is an isomorphism.  Similarly, let $A= D[X,Y; a,\der\}$ be a GWPA of rank $n\geq 1$ and $I$ be a subset of
 the set $\{ 1, \ldots , n\}$. Let $s_I$ be a bijection of the set $X\cup Y=\{ X_1, \ldots , X_1, Y_1, \ldots , Y_n\}$ which is given by the rule
 $$ s_I(X_i)= \begin{cases}
Y_i& \text{if }i\in I,\\
X_i& \text{if }i\not\in I, \\
\end{cases} \;\; {\rm and}\;\; s_I(Y_i)= \begin{cases}
X_i& \text{if }i\in I,\\
Y_i& \text{if }i\not\in I. \\
\end{cases}  $$
Let $\sign (I)\der := (\varepsilon_1 \der_1, \ldots ,\varepsilon_n\der_n )$ where $\varepsilon_i=\begin{cases}
-1& \text{if }i\in I,\\
1& \text{if }i\not\in I. \\
\end{cases} $
Then the $D$-homomorphism  of Poisson algebras
\begin{equation}\label{sIiso}
s_I: A\ra D[s_I(X), s_I(Y); a, \sign (I)\der \} , \;\; X_i\mapsto s_I(X_i), \;\; Y_i\mapsto s_I(Y_i), \;\; d\mapsto d \;\; (d\in D),
\end{equation}
is an isomorphism.

 Recall that $\d_{ij}$ is the {\em Kronecker  delta function}. The next proposition shows that the Poisson algebras $D[X,Y; \der , \alpha ]$ are GWPAs.
\begin{proposition}\label{b10Apr16}
The Poisson algebra $\CA = D[X,Y; \der , \alpha ]$ is a GWPA of rank $n$  $$ D[H_1, \ldots , H_n] [X,Y; H , \der \}$$ where $D[H_1, \ldots , H_n]$ is  a Poisson polynomial algebra over $D$  such that $\{ H_i, D\} =0$ and $\{ H_i, H_j\} =0$ for all $i,j$, $H=(H_1, \ldots , H_n)$ and  $\der_i(H_j) =\d_{ij}\alpha_jH_j$ for all $i,j$.
\end{proposition}

{\it Proof}. Consider the following elements of the polynomial algebra $\CA =D[X,Y]$, 
$$H_1=X_1Y_1, \ldots , H_n=X_nY_n.$$ Then $\{ H_i, D\}=0$ and $\{ H_i, H_j\} =0$ for all $i,j$. So,  the elements $H_1, \ldots , H_n$ belong to the absolute centre of the Poisson algebra $\CD = D[H_1, \ldots , H_n]$. Let $A= D[H_1, \ldots , H_n] [X,Y; H , \der \} $. It follows from the defining relations of the Poisson algebras $\CA$ and $A$ that there is an epimorphism $\CA \ra A$ of Poisson algebras given by the rule $X_i\mapsto X_i$, $Y_i\mapsto Y_i$, $d\mapsto d$ where $d\in D$ (since $X_iY_i\mapsto H_i$) which is clearly a bijection (it is the `identity map' of associative algebras when we identify $X_iY_i$ with $H_i$). $\Box $

$\noindent $

By Proposition \ref{b10Apr16}, the Poisson algebra $\CA = D[X,Y; \der , \alpha]=\oplus_{\beta \in \Z^n} \CA_\beta$ is $\Z^n$-graded ($\CA_\beta\CA_\g \subseteq \CA_{\beta +\g}$ and
 $\{ \CA_\beta , \CA_\g\} \subseteq \CA_{\beta+\g}$ for all $\beta , \g \in \Z^n$) where $\CA_\beta = \CD v_\beta$, $\CD =  D[H_1, \ldots , H_n]$ and $v_\beta = \prod_{i=1}^n v_{\beta_i}(i)$ where $v_j(i)=\begin{cases}
X_i^j& \text{if }j>0,\\
1& \text{if }j=1,\\
Y_i^{|j|}& \text{if }j<0.\\
\end{cases}$

{\bf Examples of GWPAs.} 1. If $D$ is an arbitrary algebra with {\em trivial Poisson bracket} (i.e. $\{ \cdot , \cdot \}=0$) then $\CZ (D) = Z(D)$ and the condition $a\in \CZ (D)^n $ in the Definition of GWPA means $a\in Z (D)^n$.  If, in addition, $D$ is a commutative algebra then $\CZ (D) = D$ and the condition $a\in \CZ (D)^n $ in the Definition of GWPA  is redundant. So, if $D$ is a {\em commutative} algebra with trivial Poisson bracket then any choice of elements $a=(a_1, \ldots , a_n)$ and $\der = (\der_1, \ldots , \der_n)\in \Der_K(D)^n$ such that $\der_i(a_j) = 0$  for all $i\neq j$ determines a GWPA $D[X,Y;a, \der\}$ of rank $n$. If, in addition, $n=1$ then there is no restriction on $a_1$ and $\der_1$.

2. The {\em classical Poisson polynomial algebra} $P_{2n} = K[X_1, \ldots , X_n , Y_1, \ldots , Y_n]$ ($\{ Y_i, X_j\} =\d_{ij}$ and $\{X_i, X_j\} = \{ X_i, Y_j\} = \{ Y_i, Y_j\} =0$ for all $i\neq j$) is a GWPA
\begin{equation}\label{P2nKH}
P_{2n}=K[H_1, \ldots , H_n][X,Y; a, \der \}
\end{equation}
where $K[H_1, \ldots , H_n]$ is a Poisson polynomial algebra with trivial Poisson bracket, $a=(H_1, \ldots , H_n)$, $\der = (\der_1, \ldots , \der_n)$ and $\der_i= \frac{\der}{\der_{H_i}}$ (via the  isomorphism of Poisson algebras $P_{2n}\ra K[H_1, \ldots , H_n][X,Y; a, \der \}$, $X_i\mapsto X_i$, $Y_i\mapsto Y_i$).

3. $A= D[X,Y; a,\der \}$ where $D=K[H_1, \ldots , H_n]$ is a Poisson polynomial algebra with trivial Poisson bracket, $a= (a_1, \ldots , a_n)\in K[H_1]\times \cdots \times K[H_n]$, $\der = (\der_1, \ldots , \der_n)$ where  $\der_i=b_i\der_{H_i}$ (where $\der_{H_i}=\frac{\der}{\der H_i}$)  and  $b_i\in K[H_i]$. In particular, $D[X,Y; (H_1, \ldots , H_n), (\der_{H_1}, \ldots , \der_{H_n})\}=P_{2n}$ is the classical Poisson polynomial algebra.

Let $S$ be a multiplicative set of $D$. Then $S^{-1}A\simeq (S^{-1} D)[X,Y; a, \der \}$ is a GWPA. In particular, for $S=\{ H^\alpha \, | \, \alpha \in \Z^n\}$ we have $K[H_1^{\pm 1}, \ldots , H_n^{\pm 1}][X,Y; a, \der \}$. In the case $n=1$, the Poisson algebra $$K[H_1^{\pm 1}][X_1, Y_1; a_1, -H_1\frac{d}{dH_1}\}$$ where $a_1\in K[H_1^{\pm 1}]$  is,  in fact, {\em isomorphic} to  a Poisson algebra   in the paper of Cho and Ho \cite{CHO-OH} which is obtained as a quantization of a certain GWA with respect to the quantum parameter $q$.  In \cite[Theorem 3.7]{CHO-OH}, a Poisson simplicity criterion  is given for this Poisson algebra. 

4.  Let $D=K[ C, H]$ be a Poisson  polynomial algebra with trivial Poisson bracket, $a\in D$ and $\der $ is a derivation of the algebra $D$. The GWPA $A= D[X, Y; a, \der \}$ of rank 1 is a generalization of some Poisson algebras that are associated with $U({\rm sl}_2)$, see the next example.

5. Let $U=U({\rm sl}_2)$ be the universal enveloping algebra of the Lie algebra
$${\rm sl_2}=K\langle X,Y, H\, | \, [H, X]=X,\;\;  [ H, Y]=-Y, \;\; [ X,Y]=2H\rangle$$  over a field $K$ of characteristic zero. The associated graded algebra ${\rm gr} (U)$ with respect to the filtration $\CF = \{ \CF_i\}_{i\in \N}$,  that is determined by the total degree of the elements $X$, $Y$ and $H$, is a Poisson polynomial algebra $K[X,Y,H]$ where
$$ \{ H, X\} = X, \;\; \{ X, Y\} =-Y\;\; {\rm and}\;\; \{ X,Y\} = 2H.$$
The element $C=XY+H^2$ belongs to the Poisson centre of the Poisson polynomial algebra ${\rm gr}(U)$. The Poisson algebra
\begin{equation}\label{grUsl2}
{\rm gr}(U)=K[C,H][X,Y;a=C-H^2, \der_H\}
\end{equation}
is a GWPA of rank 1 where $\der_H:=\frac{\der}{\der H}$.

6. Let $U$ be the universal enveloping algebra of the {\em Heisenberg Lie algebra}
\begin{eqnarray*}
 \CH_n&=&K\langle X_1, \ldots , X_n, Y_1, \ldots , Y_n, Z\, | \, [X_i, Y_j]=\d_{ij} Z,\;  [X_i, X_j]=[Y_i, Y_j]=0\; {\rm for\; all}\;\; i,j;\\
 &&   Z\; {\rm is \; a \; Poisson \; central \; element}\rangle .
\end{eqnarray*}
The associated graded algebra ${\rm gr}(U)$ with respect to the filtration by the total degree of the canonical generators is a Poisson polynomial algebra $K[X_1, \ldots , X_n, Y_1, \ldots , Y_n, Z]$  where, for all $i,j$, $$\{ X_i, Y_j\} = \d_{ij}Z, \;\; \{ X_i, X_j\} = \{ Y_i, Y_j\}=0$$ and the element $Z$ belongs to the Poisson centre of ${\rm gr}(U)$. Then the polynomial algebra
\begin{equation}\label{grUHeis}
{\rm gr}(U)= D[X,Y;a, \der\}
\end{equation}
is a GWPA of rank $n$ where $D= K[H_1, \ldots , H_n, Z]$ is a Poisson polynomial algebra with trivial Poisson bracket, $X= (X_1, \ldots , X_n)$, $Y=(Y_1, \ldots , Y_n)$, $a=(a_1=H_1, \ldots , a_n=H_n)$,  $\der = (Z\der_{H_1}, \ldots , Z\der_{H_n})$ and $\der_{H_i}:=\frac{\der}{\der H_i}$.

Let $A_s= D_s[X_{(s)},Y_{(s)}; a_{(s)}, \der_{(s)}\}$ be GWPAs of rank $n_s$ where $s=1, \ldots , m$. The tensor product of algebras
\begin{equation}\label{tpmAs}
A=\bigotimes_{s=1}^m A_s=\bigg( \bigotimes_{s=1}^m D_s\bigg) [X,Y; a, \der\}
\end{equation}
is a GWPA of rank $n_1+\cdots + n_m$ where $X=(X_{(1)}, \ldots , X_{(m)})$, $Y=(Y_{(1)}, \ldots , Y_{(m)})$, $a=(a_{(1)}, \ldots , a_{(m)})$ and $\der=(\der_{(1)}, \ldots , \der_{(m)})$. The Poisson structure on $A$ is a tensor product of Poisson structures on $A_s$, i.e. for all elements $u=\otimes_{s=1}^mu_s$, $v=\otimes_{s=1}^mv_s\in A$ (where $u_s, v_s\in A_s$),
$$ \{ u,v\} = \sum_{s=1}^mu_1v_1\t\cdots \t \{ u_s, v_s\} \t\cdots \t u_mv_m.$$
{\it Example.} The classical Poisson polynomial algebra $P_{2n}$ (see (\ref{P2nKH})) is the tensor product $P_2^{\t n}$ of $n$ copies of the classical Poisson polynomial algebra $P_2$.

The {\em opposite algebra} $A^{op}$ of an associative algebra $A$ is an algebra $A^{op}$ which is equal to   $A$ as a vector space and the product in $A^{op}$ is given by the rule $a\cdot b= ba$. If the algebra $A$ is a Poisson  algebra then so is its opposite algebra $A^{op}$ where the bracket is the same. Let $A= D[X,Y; a, \der\}$ be a GWPA of rank $n$. Then the {\em opposite Poisson algebra} to $A$,
\begin{equation}\label{AopD}
A^{op}=D^{op} [ X,Y; a, \der\},
\end{equation}
is a GWPA of rank $n$.

{\bf An algebraic torus action on a GWPA.} Let $A= D[X,Y; a, \der\}$ be a GWPA of rank $n$ and $\AutPois (A)$ be the group of automorphisms of the Poisson algebra $A$. Elements of $\AutPois (A)$ are called {\em Poisson automorphisms} of $A$.  For each element $\l = (\l_1, \ldots , \l_n)\in K^{*n}$, the $K$-algebra homomorphism
$$t_\l : A\ra A, \;\; X_i\mapsto \l_i X_i, \;\; Y_i\mapsto \l_i^{-1} Y_i, \;\; d\mapsto d\;\; (d\in D),$$ is an automorphism of the Poisson algebra $A$.  The subgroup $\mT^n = \{ t_\l \, | \, \l \in K^{*n}\}$ of $\AutPois (A)$ is an {\em algebraic torus} $\mT^n \simeq K^{*n}$, $t_\l \mapsto \l$. For all $\alpha \in \Z^n$ and  $u_\alpha \in A_\alpha = Dv_\alpha$, $t_\l (u_\alpha ) = \l^\alpha \cdot u_\alpha$ where $\l^\alpha = \prod_{i=1}^n \l_i^{\alpha_i}$.

The subgroup $$\AutPois (D)^{\der, a}:=\{ \s \in \AutPois (D) \, | \, \s \der_i = \der_i \s\;\; {\rm  and} \;\; \s (a_i) = a_i\;\; {\rm  for}\;\; i=1, \ldots , n\}$$ of $\AutPois (D)$ can be seen as a subgroup of $\AutPois (A)$ where each automorphism $\s \in \AutPois^{\der , a}(D)$  trivially acts at $X$ and $Y$, i.e. $\s (X_i) = X_i$ and $\s (Y_i) = Y_i$.  Clearly,
\begin{equation}\label{TnAut}
\mT^n\times \AutPois (D)^{\der, a} \subseteq \AutPois (A).
\end{equation}

{\bf Associated graded algebra of a GWA is a GWPA.} Let $A=D[X,Y; \s , a]$ be a GWA of rank $n$ such that $D=\cup_{i\in \N}D_i$ is a filtered algebra ($D_iD_j\subseteq D_{i+j}$ for all $i,j\in \N$; $D_{-1}=0$), 
$$[d_i, d_j]\in D_{i+j-\nu}\;\; {\rm  for \; all}\;\; d_i\in D_i\;\; {\rm  and}\;\; d_j\in D_j\;\; {\rm   where}\;\;  \nu\;\; {\rm  is\;  a \; positive \; integer;}$$
 $\s_i(D_j) = D_j$ and $(\s_i-1)(D_j) \subseteq D_{j-\nu}$ for all $i=1, \ldots , n$ and $j\in \N$. Suppose that $a_i\in D_{d_i}\backslash D_{d_i-1}$ for some $d_i\geq 1$. The algebra $A$ admits a filtration $\{ A_s\}_{s\in \frac{1}{2}\N}$ where $$A_s=\sum_{i+d\cdot \alpha \leq s}D_iv_\alpha , \;\;  d=(d_1, \ldots , d_n), \;\; \alpha = (\alpha_1, \ldots , \alpha_n)\in \Z^n\;\; {\rm  and} \;\; d\cdot \alpha = \frac{1}{2}\sum_{i=1}^n d_i|\alpha_i|.$$ The associated graded algebra
$$ {\rm gr} (A) = {\rm gr}(D)[X, Y; (\id , \ldots , \id ), \oa )={\rm gr}(D)[X_1,\ldots , X_n,  Y_1, \ldots , Y_n]/(X_1Y_1-\oa_1, \ldots , X_nY_n-\oa_n)$$
is a {\em commutative} GWA where $\oa_i = a_i+D_{d_i-1}\in D_{d_i}/D_{d_i-1}$.  For all elements $u_s\in A_s$ and $u_t\in A_t$,
\begin{equation}\label{usut}
[u_s, u_t]\in A_{s+t-\nu}.
\end{equation}
Let $\bu_s= u_s+A_{s-1}\in A_s/A_{s-1}$ and $\bu_t= u_t+A_{t-1}\in A_t/A_{t-1}$. The bracket
$$\{ u_s, u_t\} :=\overline{[u_s, u_t]}:=[u_s, u_t]+  A_{s+t-\nu -1}\in A_{s+t-\nu }/A_{s+t-\nu -1}$$
determines the Poisson structure on ${\rm gr}(A)$. For each $i=1, \ldots , n$, the map
$$\der_i:= \overline{\s_i-1}:{\rm gr} (D)\ra {\rm gr}(D),\;\; {\rm gr}(D)_j\ni \ob_j\mapsto (\s_i-1) (b_j) +D_{j-\nu -1}\in {\rm gr}(D)_{j-\nu },$$
is a $K$-derivation of the commutative algebra ${\rm gr}(D)$. The derivations $\der_1, \ldots , \der_n$ commute since the automorphisms $\s_1, \ldots , \s_n$ commute. Notice that $$ [X_i, b_j]=(\s_i-1) (b_j) X_i\;\; {\rm and}\;\; [Y_i, b_j]=(\s_i^{-1}-1) (b_j) Y_i.$$
Hence, $\{ X_i, \ob_j\} = \der_i(b_j) X_i$ and $\{ Y_i, \ob_j\} = -\der_i(b_j) Y_i$ since
$$(\s_i^{-1} -1) (b_j) = -(\s_i-1)\s_i^{-1} (b_j) \equiv -\der_i(\ob_j)\mod D_{j-\nu -1}.$$ Therefore, the Poisson algebra ${\rm gr} (A)$ is a GWPA ${\rm gr}(D)[X, Y; \oa, -\der \}$ where $ \oa= (\oa_1, \ldots , \oa_n)$ and $ -\der = (-\der_1, \ldots , -\der_n)$. So, we proved that the following proposition  holds.

\begin{proposition}\label{c10Apr16}
Let $A=D[X,Y; \s , a]$ be a GWA of rank $n$ such that $D=\cup_{i\in \N}D_i$ is a filtered algebra; $[d_i, d_j]\in D_{i+j-\nu}$ for all $d_i\in D_i$ and $d_j\in D_j$ where $\nu$ is a positive integer; $\s_i(D_j) = D_j$ and $(\s_i-1)(D_j) \subseteq D_{j-\nu}$ for all $i=1, \ldots , n$ and $j\in \N$. Suppose that $a_i\in D_{d_i}\backslash D_{d_i-1}$ for some $d_i\geq 1$. Let $\{ A_s\}_{s\in \frac{1}{2}\N}$ be the filtration as above. The associated graded algebra ${\rm gr}(A)$ is a GWPA ${\rm gr}(D)[X, Y; \oa , -\der \}$ where $\oa$ and $-\der$ are defined above.
\end{proposition}

{\it Examples.} 1. The $n$'th Weyl algebra $A_n$ is a GWA $K[H_1, \ldots , H_n][X,Y; \s , a]$ where $\s_i(H_j) = H_j-\d_{ij}$ and $a_i= H_i$ for $i,j=1, \ldots , n$. The polynomial algebra $D= K[H_1, \ldots , H_n]$ admits  a natural filtration $\{ D_i\}_{i\in \N}$ by the total degree of the variables $H_1, \ldots , H_n$. The automorphisms $\s_1, \ldots , \s_n$ satisfy the conditions of Proposition \ref{c10Apr16} with $\nu = 1$, $d_1=\cdots = d_n=1$ and $ \der_1=-\frac{\der}{\der H_1}, \ldots ,\der_n=-\frac{\der}{\der H_n}$. Notice that ${\rm gr}(D) = D$. By Proposition \ref{c10Apr16}, the algebra
$$ {\rm gr}(A_n) \simeq D[X,Y]/(X_1Y_1-H_1, \ldots , X_nY_n-H_n)\simeq K[X,Y]=P_{2n}$$ is a GWPA
$D[X,Y; (H_1, \ldots , H_n), (\frac{\der}{\der H_1}, \ldots ,\frac{\der}{\der H_n} )\}$ which is the classical Poisson algebra $P_{2n}$ with the canonical Poisson bracket ($\{ Y_i, X_j\} = \d_{ij}$,  $\{ X_i, X_j\} = \{ X_i, Y_j\} = \{ Y_i, Y_j\} =0$ for all  $i,j$ such that $i\neq j$).

2. The universal enveloping algebra $U= U({\rm sl}_2)$ is the GWA $A= K[C,H][X,Y; \s , a]$ of rank 1 where $\s (H) = H-1$, $\s (C) = C$ and $a= C-H(H+1)$ (the element $C$ is the {\em Casimir element}, $C= YX+H(H+1)$). The filtration $\CF = \{ \CF_i \}_{i\in \N}$  on $U$ that was considered above (which is defined by the total degree of the canonical generators $X$, $Y$ and $H$) induces a filtration $\{ D_i:= D\cap \CF_i \}_{i\in \N}$ on the polynomial algebra $D=K[C,H]$. Clearly,
 $$D_i = \bigoplus_{2s+t\leq i}KC^sH^t\;\; {\rm  for\;  all}\;\; i\in \N .$$ The automorphism $\s$ and the filtration $\{ D_i\}_{i\in \N}$ satisfy the conditions of Proposition \ref{c10Apr16} where $d_1=2$ and $\nu = -1$. The associated  graded Poisson algebra ${\rm gr}(A)\simeq K[C, H][X, Y; C-H^2, \der_H\}$
 is canonically isomorphic to the associated   graded Poisson algebra ${\rm gr}(U)$ as $\N$-graded  Poisson algebra (since ${\rm gr}(A)_{\frac{1}{2}+i}=0$ for all $i\in \N$), see (\ref{grUsl2}).

 The filtration $\{ D'_i:=\oplus_{j\leq i}K[C]H^i\}_{i\in \N}$ also satisfies the conditions of Proposition \ref{c10Apr16} where $d_1=2$ and $\nu = -1$ but the associated graded algebra ${\rm gr}'(A)$ is a GWPA $K[C,H][X,Y; -H^2, \der_H\}$. The associated  graded Poisson  algebras ${\rm gr}(A)$ and ${\rm gr}'(A)$ are not isomorphic  since the algebra ${\rm gr}(A)$  is smooth but the algebra  ${\rm gr}'(A)\simeq K[C]\t K[X,Y](XY-H^2)$ is singular as the points $\{ (C,H,X,Y)=(\l , 0,0,0)\, | \, \l \in K\}$ are singular. So, the Poisson algebras ${\rm gr}(A)$ and ${\rm gr}'(A)$ are also not isomorphic.


\section{Poisson simplicity criterion for generalized Weyl  Poisson   algebras}\label{SIMCRPGWA}

In this section,  for  generalized Weyl Poisson algebras, a proof of the  Poisson simplicity criterion (Theorem \ref{10Apr16}) is given, an explicit descriptions of their Poisson centre and absolute centre are obtained (Proposition \ref{a13Apr16}) and  a proof of the  criterion for the absolute centre being a field (Proposition \ref{bb13Apr16}) is given.

Let $A$ be a Poisson algebra. An ideal $I$ of the associative algebra $A$ is called a {\em Poisson ideal} if $\{ A, I\} \subseteq I$. A Poisson ideal is also called an {\em ideal of the Poisson algebra}. Suppose that $\CD$ be a set of derivations of the associative algebra $A$. Then the set $A^\CD :=\{ a\in A\, | \, \der (a) =0$ for all $\der\in \CD\}$ is a subalgebra of $A$ which is called the {\em algebra of $\CD$-constants} (or the {\em algebra of constants} for $\CD$). An ideal $J$ of the algebra $A$ is called a $\CD$-{\em invariant ideal} if $\der (J)\subseteq J$ for all $\der \in \CD$.

{\bf The Poisson centre and the absolute centre of a GWPA.} Let $A= D[X,Y; a, \der\}$ be a GWPA of rank $n$. For all elements $\l , d\in D$, $\alpha \in \Z^n$ and $i=1, \ldots , n$
\begin{equation}\label{dllv}
\{ d, \l v_\alpha\} = (-\pad_\l +\l \sum_{i=1}^n \alpha_i\der_i)(d) v_\alpha,
\end{equation}

\begin{equation}\label{dllv1}
\{ v_{\pm 1}(i), \l v_\alpha\} = \begin{cases}
\mp\der_i(\l ) v_{\alpha \pm e_i}& \text{if }\alpha_i=0 \;\; {\rm or}\;\; \sign (\alpha_i) = \pm ,\\
(\mp \der_i(\l ) a_i +\l \alpha_i\der_i(a_i) ) v_{\alpha \pm e_i}& \text{if } \sign (\alpha_i) = \mp .\\
\end{cases}
\end{equation}
The next proposition describes the centre, the Poisson centre and the absolute centre of a GWPA.

\begin{proposition}\label{a13Apr16}
Let $A= D[X,Y; a, \der\}$ be a GWPA of rank $n$. Then
\begin{enumerate}
\item $Z(A) = Z(D)[X,Y]$.
\item $\PZ (A) = \bigoplus_{\alpha \in \Z^n} \PZ (A)_\alpha$ is a $\Z^n$-graded (associative) algebra where $\PZ (A)_\alpha = D_\alpha v_\alpha$, $D_0=\PZ (D)^\der$ and, for all $\alpha \neq 0 $, $D_\alpha = \{ \l \in D^\der\, | \, \pad_\l = \l \sum_{i=1}^n \alpha_i \der_i, \; \l \alpha_i\der_i(a_i)=0$ for $i=1, \ldots , n\}$.
\item  $\CZ (A) = \bigoplus_{\alpha \in \Z^n} \CZ (A)_\alpha$ is a $\Z^n$-graded (associative) algebra where  $\CZ_\alpha = D_{[\alpha]}v_\alpha$,  $\CZ (A)_0= \CZ (D)^\der$ and, for all $\alpha \neq 0 $, $D_{[\alpha]}=Z(D)\cap D_\alpha$.
\end{enumerate}
\end{proposition}

{\it Proof}. 1. Statement 1 is obvious.

2. The GWPA $A=\oplus_{\alpha \in \Z^n} A_\alpha$ is a $\Z^n$-graded Poisson  algebra, hence so is its Poisson centre, i.e. $\PZ (A) = \oplus_{\alpha \in \Z^n} \PZ (A)_\alpha$ where $\PZ (A)_\alpha = \PZ (A)\cap A_\alpha $. Since $A_\alpha = Dv_\alpha$ for all $\alpha \in \Z^n$, statement 2 follows from (\ref{dllv}) and (\ref{dllv1}).

3. Statement 3 follows from statements 1 and 2. $\Box $


The next corollary shows that, in general, the Poisson centre and the absolute  centre of a GWPA $A$ is small.

\begin{corollary}\label{a17Apr16}
Let $A= D[X,Y;a,\der\}$ be a GWPA of rank $n$. Suppose that char$(K)=0$ and the elements $\der_1(a_1), \ldots , \der_n(a_n)$ are non-zero-divisors in the algebra $D$ (eg, $D$ is a domain and $\der_1(a_1)\neq 0, \ldots , \der_n(a_n)\neq 0$). Then
\begin{enumerate}
\item $\PZ (A) = \PZ (D)^\der$.
\item $\CZ (A) = \CZ (D)^\der$.
\end{enumerate}
\end{corollary}

For an element $\alpha = (\alpha_1, \ldots , \alpha_n)\in \Z^n$,  the set $\supp (\alpha ) :=\{ i\, | \, \alpha_i\neq 0\}$ is called the {\em support} of $\alpha$.

\begin{corollary}\label{b17Apr16}
Let $A= D[X,Y;a,\der\}$ be a GWPA of rank $n$. Suppose that char$(K)=0$. Then, for all elements $\alpha \in \Z^n\backslash \{ 0\}$, $D_\alpha \subseteq D^{\der , \pad (\der (a))}\cap \ann_D \{\der_i(a_i)\, | \, $  $i\in \supp (\alpha )\}$, i.e.
\begin{enumerate}
\item $\{ D_\alpha , a_i\} =0$ for $i=1, \ldots , n$, and
\item $D_\alpha \der_i (a_i)=0$ for all $i$ such that $\alpha_i\neq 0$.
\end{enumerate}
\end{corollary}

{\it Proof}. By Proposition \ref{a13Apr16}.(3), $D_\alpha \der_i (a_i)=0$ for all $i\in \supp (\alpha )$ (since char$(K)=0$). Then, for all $\l \in D_\alpha$ and $i=1, \ldots , n$, $\{ \l , a_i\} = \pad_\l (a_i) = \sum_{i=1}^n \l \alpha_i \der_i(a_i) =0$, i.e. $\{ D_\alpha , a_i\} =0$ for $i=1, \ldots , n$. $\Box$

Let $A=\bigoplus_{i\in \Z} A_i$ be a $\Z$-graded (associative) algebra. Each element $a\in A$ is a unique sum $a=\sum_{i\in \Z} a_i$ where $a_i\in A_i$. The {\em length} $l(a)$ of the element $a$ is equal to $-\infty$ if $a=0$, and, for $a\neq 0$,  $l(a) := n-m$ where $n=\max \{ i \, | \, a_i\neq 0\}$ and $m=\min \{ i \, | \, a_i\neq 0\}$.

Let $A$ be a Poisson algebra and $z\in \CZ (A)$. The $zA$ is a Poisson ideal of $A$. {\em If the Poisson algebra $A$ is simple then necessarily the absolute centre $\CZ (A)$ is a field.}\\

{\bf Proof of Proposition \ref{bb13Apr16}}. $(\Rightarrow )$ Suppose that $p={\rm char}(K)\neq 0$. Then, by Proposition \ref{a13Apr16}, the element $1+X^p$ of $\CZ (A)$ is not invertible. Therefore, we must have $p=0$. The algebras $A$ and $\CZ (A)$ are $\Z^\alpha$-graded algebras and $\CZ (A)_0= \CZ (D)^\der$. Therefore, $\CZ (D)^\der$ must be a field.

Suppose that $D_{[\alpha ]}\neq 0$ for some $\alpha \neq 0$. Then $\alpha_i\neq 0$  for some $i$. Fix a nonzero element of $\CZ(A)_\alpha =D_{[\alpha ]}v_\alpha$, say $\l v_\alpha$ where $\l \in D_{[\alpha ]}$. Since $\l v_\alpha$ is a unit, $(\l v_\alpha )^{-1}=\mu v_{-\alpha }$ (since the algebra $A$ is a $\Z^n$-graded algebra), and so 
$$1= \l v_\alpha \cdot \mu v_{-\alpha}=\l \mu a^{|\alpha |}\;\; {\rm  and}\;\;  1=\mu v_{-\alpha}\cdot \l v_\alpha = \mu \l a^{|\alpha |}$$ 
where $ a^{|\alpha |}:=\prod_{i=1}^na_i^{|\alpha_i|}\in \CZ (A)$. Hence, $a^{|\alpha |}$ is a unit in $\CZ (A)$, then the elements $\l $ and $\mu$ are units in $D$. Clearly, $v:=1+\l v_\alpha \in \CZ (A)$. The algebra $A$ is a $\Z^n$-graded algebra. In particular, it is a $\Z e_i$-graded algebra (since $\Z e_i \subseteq \Z^n$). Let $l_i$ be the length with respect to the $\Z e_i$-grading (which is a $\Z$-grading). Then, for all nonzero elements $u\in A$,
$$ l_i(uv) = l_i(vu) = l_i(u) +l_i(v) \geq l_i(v)=|\alpha_i|>0,$$
since the elements 1 and $\l$ are units.
This implies that the element $u$ is not a unit. Therefore, $D_{[\alpha ]}=0$ for all $\alpha \in \Z^n \backslash \{ 0\}$, by Proposition \ref{a13Apr16}.(3).

$(\Leftarrow )$ By Proposition \ref{a13Apr16}.(3), $\CZ (A) = \CZ (D)^\der$ is a field.  $\Box $


An ideal $I$ of an algebra $A$ is called a {\em proper ideal} if $I\neq 0 , A$.\\

{\bf  Proof of Theorem \ref{10Apr16}}.  $(\Rightarrow )$ Suppose that $\ga$ is a proper $\der$-invariant Poisson ideal of the Poisson algebra $D$ then $\ga A=\oplus_{\alpha \in \Z^n}\ga v_\alpha$ is a proper ideal of the Poisson algebra $A$. So, the first condition holds.

Suppose that $\gb := Da_i+D\der_i(a_i) \neq  D$ for some $i$. Then 
$$I=\bigoplus_{\alpha\in \Z^n, \alpha_i\neq 0}Dv_\alpha \oplus \bigoplus_{\alpha \in \Z^n, \alpha_i = 0} \gb  v_\alpha$$ is a proper ideal of the Poisson algebra $A$. So, the second condition holds.

The third condition obviously holds (if a nonzero element $z$ of $\CZ (A)$ is also a non-unit then $zA$ is a proper Poisson ideal of $A$).

 $(\Leftarrow )$  Suppose that conditions 1 and 2 hold. Then the implication follows from the Claim.

 {\em  Claim. Suppose that conditions 1 and 2 hold. Then every nonzero Poisson ideal of $A$ intersects nontrivially $\CZ (A)$.}

Let $I$ be a nonzero Poisson ideal $A$. We have to show that $I\cap \CZ (A) \neq 0$. Let $u=\sum_{\alpha \in \Z^n} u_\alpha$ be a nonzero element of $I$ where $u_\alpha \in A_\alpha$. The set $\supp (u) = \{ \alpha \in \Z^n\, | \, u_\alpha \neq 0\}$ is called the {\em support} of $u$. Recall that, for $\alpha \in \Z^n$, $|\alpha |=\alpha_1+\cdots +\alpha_n$. The additive group $\Z^n$ admits the {\em degree-by-lexicographic ordering} $\leq $ where $\alpha <\beta$ iff either $|\alpha | <|\beta |$ or $|\alpha | =|\beta |$ and there exists an element $i\in \{ 1, \ldots , n\}$ such that $\alpha_j= \beta_j$ for all $j<i$ and $\alpha_i<\beta_i$. Clearly, the inequalities $\alpha \leq \beta$ and $\beta \leq \alpha$ are equivalent to the equality  $\alpha = \beta$. The partially ordered set $(\Z^n, \leq )$ is a {\em linearly ordered} set (for all distinct elements $\alpha , \beta \in \Z^n$ either $\alpha >\beta$ or $\alpha <\beta$) and $\alpha < \beta$ implies that $\alpha +\g <\beta +\g$ for all $\g \in \Z^n$. Every nonzero element $b=\sum_{\alpha \in \Z^n} b_\alpha$ of $A$ (where $b_\alpha \in A_\alpha $) can be written as $$b= b_\alpha +\cdots$$ where $\alpha$ is the maximal element of $\supp (b)$ and the three dots denote smaller terms (i.e. the sum $\sum_{\beta <\alpha }b_{\beta}$). The term $b_\alpha = \l_\alpha v_\alpha$ is called the {\em leading term} of $b$, denoted $\lt (b)$, and the element $\l_\alpha \in D$ is called the {\em leading coefficient} of $b$, denoted $\lc (b)$.  Since the algebra $A$ is a $\Z^n$-graded Poisson algebra,  for all nonzero elements $b,c\in A$,
\begin{equation}\label{ltbc}
\lt (bc) = \lt (b) \lt (c)
\end{equation}
provided $\lc (b) \lc (c)\neq 0$, and
\begin{equation}\label{ltbc1}
\lt (\{ b,c\} ) = \{ \lt (b),  \lt (c)\}
\end{equation}
provided $\{ \lt (b),  \lt (c)\} \neq 0$.

Up to isomorphism in (\ref{s1iso}) (i.e. interchanging some $X_i$ and $Y_i$, if necessary), we can assume that the ideal $I$ contains a nonzero element $u=\l _\alpha X^\alpha +\cdots$ where $\alpha_1\geq 0, \ldots , \alpha_n\geq 0$. Then the set of leading coefficients
$$ \ga = \{ \l_\alpha \, | \, u = \l_\alpha X^\alpha +\cdots \in I, \;\ {\rm all}\;\; \alpha_i\geq 0\}$$ of elements of $I$ is a $\der$-invariant ideal of the ring $D$ since
\begin{eqnarray*}
d_1ud_2&= & d_1\l d_2u_\alpha +\cdots\;\;\;\;\; \;\;\; \; {\rm if}\; d_1\l d_2\neq 0\;\;  (d_1, d_2\in D),\\
  uX^\beta &=& \l_\alpha X^{\alpha +\beta}+\cdots , \\
\{ u, X_i\} &=&\der_i(\l_\alpha ) X^{\alpha +e_i}+\cdots \;\; {\rm if}\;\; \der_i(\l_\alpha )\neq 0.
\end{eqnarray*}
Therefore, by condition 1, there exists an element $u=X^\alpha +\cdots \in I$ (i.e. $\l_\alpha =1$). Then using the equalities 
$$Y_iX_i^\alpha = a_iX^{\alpha -e_i}\;\; {\rm  and}\;\;\{ Y_i, X^\alpha\} =\alpha_i\der_i (a_i) X^{\alpha_i-e_i},$$    condition 2 and the fact that char$(K)=0$ (condition 3), we can assume that $u=1+\cdots \in I$, i.e. $u=1+\sum_{\alpha <0}u_\alpha$. For a finite set $S$, we denote by $|S|$ the number of its elements.  Let
$$m = \min \{ |\supp (u)|\, | \, u=1+\cdots \in I\}.$$
We can assume that $|\supp (u)|=m$. The Poisson algebra $A$ is a $\Z^n$-graded Poisson algebra. Hence, by the choice of $m$, $0=[d,u]=\sum_{\alpha <0}[d, u_\alpha]$ for all elements $d\in D$, i.e. all $u_\alpha \in Z(A)$, and so $u\in Z(A)$. Similarly, for all elements  $d\in D$ and $i=1, \ldots , n$
\begin{eqnarray*}
0&=& \{ d,u\} =\sum_{\alpha <0} \{ d, u_\alpha\} ,\;\; {\rm i.e.}\;\; \{ d,u_\alpha\}=0, \\
0&=&  \{ X_i,u\} =\sum_{\alpha <0} \{ X_i, u_\alpha\} ,\;\; {\rm i.e.}\;\; \{ X_i,u_\alpha\}=0, \\
 0&=& \{ Y_i,u\} =\sum_{\alpha <0} \{ Y_i, u_\alpha\} ,\;\; {\rm i.e.}\;\; \{ X_i,u_\alpha\}=0,
\end{eqnarray*}
i.e. all $u_\alpha \in \PZ (A)$, and so $u\in \PZ (A)$. Then $0\neq u \in Z(A)\cap \PZ (A)=\CZ (A)$, as required. $\Box$

\begin{corollary}\label{a16Apr16}
Let $A=D[X,Y;a, \der\}$ be a GWPA of rank $n$. Suppose that the conditions 1 and 2 of Theorem \ref{10Apr16} hold. Then every nonzero Poisson ideal of $A$ intersects $\CZ (A)$ nontrivially.
\end{corollary}

{\it Proof}. The corollary is precisely the Claim in the proof of Theorem \ref{10Apr16}.
$\Box $

\begin{corollary}\label{a19Apr16}
Let $D=K[H_1, \ldots , H_n]$ be a Poisson polynomial algebra with trivial Poisson bracket, $a=(a_1, \ldots , a_n)$ where $a_i\in K[H_i]$ and $\der = (b_1\der_{H_1}, \ldots , b_n\der_{H_n})$ where $b_i\in K[H_i]$. Then  the GWPA   $A=D[X,Y;a, \der\}$ of rank $n$   is a simple Poisson algebra iff char$(K)=0$, $b_1, \ldots , b_n\in K^*:=K\backslash \{ 0\}$ and $K[H_i]a_i+K[H_i]\frac{da_i}{dH_1}=K[H_i]$ for $i=1, \ldots , n$.
\end{corollary}

{\it Proof}. The corollary follows from  Theorem \ref{10Apr16}. In more detail, condition 2 of Theorem \ref{10Apr16} is equivalent to the conditions  $K[H_i]a_i+K[H_i]\frac{da_i}{dH_1}=K[H_i]$ for $i=1, \ldots , n$ (since $a_i\in K[H_i]$). Condition 1 of Theorem \ref{10Apr16} is equivalent to the condition char$(K)=0$ and  $b_1, \ldots , b_n\in K^*:=K\backslash \{ 0\}$ (since $b_iD$ is a $\der$-invariant ideal of $D$). If conditions 1 and 2 hold then condition 3 of Theorem \ref{10Apr16} holds automatically since $D^\der =K=\CZ (D)$ (then $D_{[\alpha]}=0$ for all $\alpha\in \Z^n\backslash \{ 0\}$). $\Box $

By Corollary \ref{a19Apr16}, the classical Poisson polynomial algebra $$P_{2n}\simeq K[H_1, \ldots , H_n][X,Y;(H_1, \ldots , H_n), (\der_{H_1}, \ldots , \der_{H_n})\}$$ is a simple Poisson algebra.



\small{

Department of Pure Mathematics

University of Sheffield

Hicks Building

Sheffield S3 7RH

UK

email: v.bavula@sheffield.ac.uk}

\end{document}